\documentclass[11pt,oneside,a4paper]{article}
\usepackage[T1]{fontenc}
\usepackage{lmodern}
\usepackage{xparse}
\usepackage{amsmath}
\usepackage{amssymb}
\usepackage{comment}
\usepackage[english]{babel}
\usepackage[utf8]{inputenc}
\usepackage[]{nicefrac}
\usepackage{amsmath, amssymb, amsthm}
\usepackage{graphicx}
\usepackage{tikz}
\usepackage[affil-it]{authblk}
\usepackage[]{xcolor}
\usepackage{leftidx}
\usepackage[numbers]{natbib}
\usepackage{pgfplots}
\usepackage[linesnumbered,ruled,vlined]{algorithm2e}
\usepackage[]{hyperref}
\usepackage{xcolor}
\hypersetup{%
    colorlinks = true,
    allcolors = gray,
}

\newcommand{\CG}[1]{\textcolor{blue}{#1}}
\definecolor{officegreen}{rgb}{0.0, 0.5, 0.0}

\newcommand{\dVblank}{\, {\rm d}V}

\newcommand{\mcF}{\mathcal{F}}
\newcommand{\mcE}{\mathcal{E}}

\newcommand{\mcX}{\mathcal{X}}

\newcommand{\mcL}{\mathcal{L}}

\newcommand{\mcA}{\mathcal{A}}

\newcommand{\mbR}{\mathbb{R}}
\newcommand{\mbRd}{{\mathbb{R}^d}}

\newcommand{\Omg}{{\Omega}}
\newcommand{\Gam}{{\Gamma}}
\newcommand{\gam}{{\gamma}}

\def \ub{\mathbf{u}}

\def \xb{{\bf x}}
\def \yb{{\bf y}}

\NewDocumentCommand \dV{ o }{
    \IfNoValueTF{#1}{\dVblank}
    {
        \dV_{#1}
    }
}

\NewDocumentCommand \mc{ m }{
    \mathcal{#1}
}

\NewDocumentCommand \eref{ m }{
    (\ref{eqn:#1})
}

\NewDocumentCommand \an{ m }{
    \langle {#1} \rangle
}

\NewDocumentCommand \mbf{ m }{
    \mathbf{#1}
}

\NewDocumentCommand \bs{ m }{
    \boldsymbol{#1}
}

\NewDocumentCommand \stateOne{ m }{
    \underline{\mbf{#1}}
}

\NewDocumentCommand \stateTwo{ m o }{
    \IfNoValueTF{#2}{\stateOne{#1}}
    {
    \stateOne{#1}\an{#2}
    }
}

\NewDocumentCommand \s{ m o o }{
    \IfNoValueTF{#3}{\stateTwo{#1}[#2]}
    {
    \underline{\mbf{#1}}(\mbf{#3})\an{#2}
    }
}

\NewDocumentCommand \stateOneI{ m m }{
    \underline{#1}_{#2}
}

\NewDocumentCommand \stateTwoI{ m m o }{
    \IfNoValueTF{#3}{\stateOneI{#1}{#2}}
    {
    \stateOneI{#1}{#2}\an{#3}
    }
}

\NewDocumentCommand \si{ m m o o }{
    \IfNoValueTF{#4}{\stateTwoI{#1}{#2}[#3]}
    {
        \stateOneI{#1}{#2}(\mbf{#4})\an{#3}
    }
}

\title{A FETI approach to domain decomposition for meshfree discretizations of nonlocal problems}
\date{\today}

\author{Xiao Xu\thanks{The Oden Institute for Computational Engineering and Sciences, The University of Texas at Austin, TX, xiaoxu42@utexas.edu}
\and Christian Glusa \thanks{Computer Science Research Institute, Sandia National Laboratories, NM, caglusa@sandia.gov}
\and Marta D'Elia \thanks{Computational Science and Analysis, Sandia National Laboratories, CA, mdelia@sandia.gov}
\and John T. Foster\thanks{The Oden Institute for Computational Engineering and Sciences, The University of Texas at Austin, TX, john.foster@utexas.edu}}

\begin{document}
\maketitle
\begin{abstract}
We propose a domain decomposition method for the efficient simulation of nonlocal problems. Our approach is based on a multi-domain formulation of a nonlocal diffusion problem where the subdomains share ``nonlocal'' interfaces of the size of the nonlocal horizon. This system of nonlocal equations is first rewritten in terms of minimization of a nonlocal energy, then discretized with a meshfree approximation and finally solved via a Lagrange multiplier approach in a way that resembles the finite element tearing and interconnect method. Specifically, we propose a distributed projected gradient algorithm for the solution of the Lagrange multiplier system, whose unknowns determine the nonlocal interface conditions between subdomains. Several two-dimensional numerical tests illustrate the strong and weak scalability of our algorithm, which outperforms the standard approach to the distributed numerical solution of the problem. This work is the first rigorous numerical study in a two-dimensional multi-domain setting for nonlocal operators with finite horizon and, as such, it is a fundamental step towards increasing the use of nonlocal models in large scale simulations. 
\end{abstract}
\section{Introduction}
Nonlocal models have become alternatives to classical, differential equation models to describe phenomena where long-range interactions at small scales affect the global behavior of the solution. In particular, they are the models of choice in several scientific and engineering applications where multiscale effects, anomalous behavior, and discontinuities in the solutions must be taken into account for reliable predictions. Applications of interest include, fracture mechanics \citep{silling2000reformulation,ha2011characteristics}, subsurface flow \citep{Benson2000,Schumer2003}, image processing \citep{DElia2019imaging,Gilboa2007}, stochastic processes \citep{DElia2017,Meerschaert2012}, and many others. 

The multiscale nature of nonlocal models manifests itself in the form of integral operators that embed length-scales in their definition. In their simplest form, they can be defined as 
\begin{equation*}
\mcL u(\xb)= 2 \int_{B_\delta(\xb)} (u(\yb)-u(\xb))\gam(\xb,\yb)\,d\yb,
\end{equation*}
where $B_\delta(\xb)$ is a norm-induced ball (e.g. Euclidean) centered at $\xb$ of radius $\delta$, also known as {\it horizon}, that determines the extent of the nonlocal, long-range interactions, described by the \emph{kernel} $\gamma$. The integral form also allows for discontinuities in the solutions, which means that regularity requirements are much less strict than for partial differential equations (PDEs) \citep{du12}.

However, despite their modelling generality, the use of nonlocal models is still hindered by several mathematical and computational challenges.  In this paper we focus on the prohibitively expensive computational cost, particularly apparent when the ratio between horizon and discretization length becomes large.  Moreover, we explore the treatment of virtual nonlocal interfaces, which, for nonlocal problems, is still in its infancy \citep{Alali2015,Capodaglio2020}. Specifically, we propose a new substructuring-based domain decomposition (DD) approach for meshfree discretizations of nonlocal equations that resembles non-overlapping substructuring DD approaches for PDEs. 

The first approach to nonlocal DD can be found in \cite{aksoylu2011variational}; here, the authors consider a substructuring scheme based on the solution of a Schur complement problem in a two-domain configuration. As such, this work, while rigorous, does not address issues such as scalability, performance, and numerical treatment of floating sub-domains. In \cite{Capodaglio2021} the authors introduce a rigorous multi-domain DD formulation, but they do not pursue a specific algorithm for its numerical solution nor do they present a numerical study. 

In this paper, we consider a meshfree discretization of the multi-domain formulation introduced in \cite{Capodaglio2021} and propose to solve the multi-domain system via a Lagrange multiplier approach, in a way that resembles the finite element tearing and interconnect (FETI) method \citep{farhat1991method}. Our main contributions are summarized below.
\begin{itemize}
    \item We introduce a meshfree discretization of the multi-domain DD formulation described in \cite{Capodaglio2021} and propose a solution method based on a Lagrange-multiplier approach.
    \item Following FETI's literature, we design an iterative solver for the solution of Lagrange-multiplier system that asymptotically converges in a constant number of iterations.
    \item We show the improved performance of nonlocal FETI compared to standard parallel solvers for the single-domain problem.
    \item We perform several weak and strong scalability studies that illustrate the efficiency of the proposed approach.
\end{itemize}
We highlight the fact that this is the first rigorous DD numerical study in a two-dimensional multi-domain setting for nonlocal operators with finite nonlocal interactions. As such, this work is the first step towards increasing the use of finite-horizon nonlocal models in large-scale simulations.

Furthermore, the impact of this work extends beyond nonlocal models and finds applicability in meshfree discretizations of PDEs. In fact, when the latter are characterized by long stencils \cite{bessa2013}, their discretization matrices feature a wide bandwidth, proportional to the number of degrees of freedom involved in the discretization of a certain point in the domain. This situation resembles nonlocal interactions, that result in discretization matrices whose bandwidth is proportional to the radius of interactions. The approach proposed in this work extends to long-stencil meshfree PDE discretizations in a straightforward manner and its implementation is part of our current work.

\paragraph{Outline of the paper} In Section \ref{sec:notation} we introduce the notation used throughout the paper and review important aspects of the multi-domain formulation introduced in \cite{Capodaglio2021}. In Section \ref{sec:FETI} we first introduce the FETI formulation and its meshfree discretization and then provide implementation details regarding the parallel projected gradient algorithm used in our computational tests. In Section \ref{sec:numer-exper} we illustrate the performance of the proposed algorithm with several two-dimensional tests. In particular, we show strong and weak scaling results and compare the performance of our algorithm with a standard parallel solver for the single-domain problem. Finally, in Section \ref{sec:conclusion} we summarize our contributions and propose future research guidelines. 

\section{Notation and previous work}\label{sec:notation}

Let $\gamma(\xb,\yb)$ be an integrable symmetric nonnegative\footnote{See \cite{Mengesha2013,xu2020deriving,Xu2020learning,You2020aaai,you2020data} for discussions on nonpositive kernels and \cite{DElia2017,Felsinger2015} for discussions on nonsymmetric kernels.} kernel, i.e., $\gamma(\xb,\yb) = \gamma(\yb,\xb)$ for $\xb,\,\yb\in\mathbb{R}^d$, with bounded support in the norm-induced ball of radius $\delta$, i.e. $\text{supp}(\gamma(\xb,\cdot))=\{\yb:|\xb-\yb|_{\ell^p}\leq\delta\}$. We refer to $\delta>0$ as the \emph{horizon} and let $p\in[1,\infty]$.
Let $\Omega\in\mbRd$ be an open bounded domain.
We define its associated {\it interaction domain} as the set of points outside $\Omega$ that interact with points inside $\Omega$, i.e. 
{\begin{equation}\label{eq:int_dom}
\Gamma= \{\yb \in  \mbRd \setminus\Omega\,:\, \exists \;\xb \in \Omega\;\; \hbox{such that} \;\; |\xb - \yb|_{\ell^p} \leq \delta \; \}.
\end{equation}}
The \emph{nonlocal Laplacian} is then given by
\begin{equation}\label{eq:laplacian}
\mcL u(\xb)= 2 \int_{\Omega\cup\Gam} (u(\yb)-u(\xb))\gam(\xb,\yb)\,d\yb.
\end{equation}
We define the strong form of a nonlocal volume-constrained Poisson problem  \citep{acta20,du12} as follows: given $f:\Omg\to\mbR$ and $g:\Gamma\to\mbR$, find $u:\Omega\cup\Gam\to\mbR$ such that
\begin{equation}\label{eq:strong_nonlocal}
\begin{cases}
\displaystyle  
-\mcL u(\xb)= f(\xb) & \xb\in\Omega, \\[2mm]
u(\xb) = g(\xb) &  \xb\in\Gamma.
\end{cases}
\end{equation}
Here, the second equation is the nonlocal counterpart of a Dirichlet boundary condition for PDEs and, as such, it is referred to as {\it Dirichlet volume constraint}. In this work, without loss of generality, we only consider Dirichlet constraints as the extension to Neumann\footnote{For definition and analysis of Neumann volume constraints we refer to \cite{Du2012} and for its numerical treatment we refer to, e.g., \cite{DEliaNeumann2019}.} constraints does not affect the proposed algorithm.

With the purpose of introducing a weak form for problem \eqref{eq:strong_nonlocal}, we define the function spaces
\begin{align}\label{globalspace}
W&= \{w \in L^2(\Omega\cup\Gamma) : |||w||| < \infty\},
& W^0 &= \{w\in W : w|_{\Gamma}=0 \},
\end{align}
where 
$$
|||w|||^2 =  \int_{\Omega\cup\Gamma}\int_{\Omega\cup\Gamma} |w(\yb) - w(\xb)|^2\gamma(\xb,\yb)d\yb  d\xb + \|w\|^2_{L^2(\Omega\cup\Gamma)}.
$$
For $u,v\in W$, we define the bilinear form $\mathcal{A}(\cdot,\cdot)$ and linear functional $\mcF(\cdot)$ as
\begin{align}\label{eq:bil1}
\mathcal{A}(u,v) &=
\int_{\Omega\cup\Gamma}  \int_{\Omega\cup\Gamma}\big(u(\yb) - u(\xb)\big)\big(v(\yb) - v(\xb)\big)\gamma(\xb,\yb)d\yb  d\xb 
\\
\mcF(v) &= \int_\Omega v(\xb)f(\xb)d\xb.
\end{align}
Then, the weak formulation of \eqref{eq:strong_nonlocal} is given by: given $f\in W'$, $g\in W_\Gamma$, and a kernel $\gamma$, find $u \in W$ such that
\begin{equation}\label{eq:weak_nonloc}
\mcA(u,v) = \mcF(v) \quad \forall\, v \in W^0
\qquad\mbox{subject to $u|_{\Gamma}=g$},
\end{equation}
where $W'$ denotes the dual space of bounded linear functionals on $W^0$ with respect to the standard $L^2$ duality pairing and $W_\Gamma$ denotes the nonlocal trace space defined as $W_\Gamma=\{w|_\Gamma\,\,:\,\, w\in W\}$. For more details, including well posedness of \eqref{eq:weak_nonloc}, see, e.g., \cite{acta20,DElia2020unified,du12}.

\subsection{A substructuring-based domain-decomposition method for nonlocal problems}

We briefly review the multi-domain formulation presented in \cite{Capodaglio2021}, only highlighting the aspects that are relevant for our work and keeping the same notation for consistency.

First, we illustrate how to decompose the domain $\Omega$, taking into account the presence of nonlocal interactions. We start from a non-overlapping, covering subdivision of $\Omega$ into $N$ subdomains $\{\widetilde\Omega_n\}_{n=1}^{N}$, as illustrated in Figure~\ref{fig:sketch}, left, and then add the overlaps necessary for the nonlocal interactions. For each subdomain we define the subdomain 
\begin{equation}\label{subdot}
   \Omega_n =
   \left\{ \xb\in\widetilde\Omega_n\,\,\,
   :\,\,\, |\yb-\xb|>\frac\delta2 \quad\forall \,\yb\in
 \Omega\setminus\widetilde\Omega_n  
\right\};\end{equation}
as illustrated in Figure~\ref{fig:sketch}, right, for $p=\infty$. We also subdivide the interaction domain $\Gamma$ accordingly, into a set of overlapping, covering subdomains 
\begin{equation}\label{subdg}
    \Gamma_n = \big\{ \xb\in\Gamma\,
    : \, \exists\;\yb\in \Omega_n \;\; \hbox{such that}
    \;\;|\yb-\xb| \le \delta 
    \;\big\}\qquad\mbox{for $n=1,\ldots,N$}
\end{equation}
and introduce the set of interaction domains
\begin{equation}\label{subdgt}
    \widehat\Gamma_n = \big\{ \xb\in \Omega\setminus\Omega_n\,
    : \, \exists\;\yb\in \Omega_n \;\; \hbox{such that}
    \;\;|\yb-\xb| \le \delta 
    \;\big\}\qquad\mbox{for $n=1,\ldots,N$}.
\end{equation}
As a result, for each $n$, $\Gamma_n\cup \widehat\Gamma_n$ consists of all the strips of thickness $\delta$ that surround $\Omega_n$, including in some instances a portion of $\Gamma$. This construction results in the following overlapping decomposition of $\Omega\cup\Gamma$:
\begin{equation}\label{overdd}
\Omega\cup\Gamma = \cup_{n=1}^{N} \Omega_n\cup\widehat\Gamma_n\cup\Gamma_n.
\end{equation}

As it is common for DD methods, we divide the set of subdomains $\{\Omega_n\cup\widehat\Gamma_n\cup\Gamma\}_{n=1}^{N}$ into two classes: {\it floating} subdomains if $\Gamma_n=\emptyset$ and {\it non-floating} subdomains if $\Gamma_n \ne \emptyset$. In the former case, the domain is endowed with a purely Neumann nonlocal volume constraint so that its associated nonlocal problem is singular, i.e. has a non-trivial nullspace. In practical contexts, the set of floating subdomains is always non-empty, making the treatment of the associated singular problems critical.
\definecolor{mypeach}{RGB}{255,229,180}
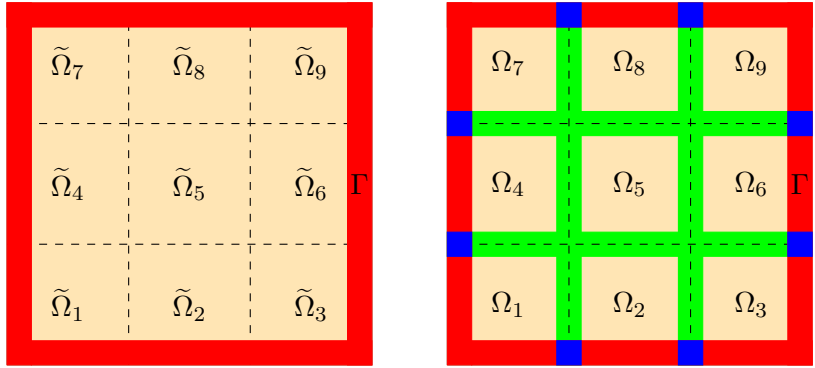
\begin{figure}
\centering
\begin{tikzpicture}[scale=0.4]
    \draw (-6,-6) -- (6,-6) -- (6,6) -- (-6,6) -- (-6,-6);
    \filldraw[fill=mypeach, draw=mypeach] (-6,-6) rectangle (6,6);
    \draw [dashed] (-2,-6) -- (-2,6);
    \draw [dashed] (2,-6) -- (2,6);
    \draw [dashed] (-6,-2) -- (6,-2);
    \draw [dashed] (-6,2) -- (6,2);
    \filldraw[fill=red, draw=red] (-6,-6) rectangle (-6+0.4*2,6);
    \filldraw[fill=red, draw=red] (-6+0.4*2,-6) rectangle (6,-6+0.4*2);
    
    \filldraw[fill=red, draw=red] (6,-6) rectangle (6-0.4*2,6);
    \filldraw[fill=red, draw=red] (6,6) rectangle (-6,6-0.4*2);
    
    \draw (6-0.4,0) node () {$\Gamma$};
    
    \draw (-4,-4) node () {$\widetilde\Omega_1$} ;
    \draw (0,-4) node () {$\widetilde\Omega_2$} ;
    \draw (4,-4) node () {$\widetilde\Omega_3$} ;
    \draw (-4,0) node () {$\widetilde\Omega_4$} ;
    \draw (0,0) node () {$\widetilde\Omega_5$} ;
    \draw (4,0) node () {$\widetilde\Omega_6$} ;
    \draw (-4,4) node () {$\widetilde\Omega_7$} ;
    \draw (0,4) node () {$\widetilde\Omega_8$} ;
    \draw (4,4) node () {$\widetilde\Omega_9$} ;
\end{tikzpicture}
\qquad
\begin{tikzpicture}[scale=0.4]
    \draw (-6,-6) -- (6,-6) -- (6,6) -- (-6,6) -- (-6,-6);
    \filldraw[fill=mypeach, draw=mypeach] (-6,-6) rectangle (6,6);
    \filldraw[fill=green, draw=green] (-2-0.4,-6) rectangle (-2+0.4,6);
    \filldraw[fill=green, draw=green] (2-0.4,-6) rectangle (2+0.4,6);
    \filldraw[fill=green, draw=green] (-6,-2-0.4) rectangle (6,-2+0.4);
    \filldraw[fill=green, draw=green] (-6,2-0.4) rectangle (6,2+0.4);
    \draw [dashed] (-2,-6) -- (-2,6);
    \draw [dashed] (2,-6) -- (2,6);
    \draw [dashed] (-6,-2) -- (6,-2);
    \draw [dashed] (-6,2) -- (6,2);
    \filldraw[fill=red, draw=red] (-6,-6) rectangle (-6+0.4*2,6);
    \filldraw[fill=red, draw=red] (-6+0.2*4,-6) rectangle (6,-6+0.4*2);
    \filldraw[fill=red, draw=red] (6,-6) rectangle (6-0.4*2,6);
    \filldraw[fill=red, draw=red] (6,6) rectangle (-6,6-0.4*2);
    
    \filldraw[fill=blue, draw=blue] (-6,2-0.4) rectangle (-6+0.4*2,2+0.4);
    \filldraw[fill=blue, draw=blue] (-6,-2-0.4) rectangle (-6+0.4*2,-2+0.4);
    \filldraw[fill=blue, draw=blue] (-2-0.4,-6) rectangle (-2+0.4,-6+0.2*4);
    \filldraw[fill=blue, draw=blue] (2-0.4,-6) rectangle (2+0.4,-6+0.4*2);
    
    \filldraw[fill=blue, draw=blue] (6,2-0.4) rectangle (6-0.4*2,2+0.4);
    \filldraw[fill=blue, draw=blue] (6,-2-0.4) rectangle (6-0.4*2,-2+0.4);
    \filldraw[fill=blue, draw=blue] (-2-0.4,6) rectangle (-2+0.4,6-0.2*4);
    \filldraw[fill=blue, draw=blue] (2-0.4,6) rectangle (2+0.4,6-0.4*2);

    \draw (6-0.4,0) node () {$\Gamma$};
    
    \draw (-4,-4) node () {$\Omega_1$} ;
    \draw (0,-4) node () {$\Omega_2$} ;
    \draw (4,-4) node () {$\Omega_3$} ;
    \draw (-4,0) node () {$\Omega_4$} ;
    \draw (0,0) node () {$\Omega_5$} ;
    \draw (4,0) node () {$\Omega_6$} ;
    \draw (-4,4) node () {$\Omega_7$} ;
    \draw (0,4) node () {$\Omega_8$} ;
    \draw (4,4) node () {$\Omega_9$} ;
\end{tikzpicture}
\caption{On the left, a non-overlapping covering subdivision of $\Omega$ (yellow) and the interaction domain $\Gamma$ (red). On the right, the associated decomposition that accounts for nonlocal interactions. Green regions correspond to $\bigcup_{n=1}^N \widehat\Gamma_n$ and blue regions highlight portions of $\Gamma$ that are shared by two subdomains.}
\label{fig:sketch}
\end{figure}

Based on the decomposition above, we define the multi-domain system; this is composed of $N$ equations defined on $\Omega_n\cup\widehat\Gamma_n\cup\Gamma_n$, $n=1,\ldots,N$. In order to deal with the fact that the domains above are overlapping, we introduce the following generalized indicator functions
\begin{equation}\label{eq:zeta}
\zeta_{\mcA}(\xb,\yb)=
\sum_{n=1}^{N}\mcX_{\Omega_n\cup\widehat\Gamma_n\cup\Gamma_n}(\xb)
\mcX_{\Omega_n\cup\widehat\Gamma_n\cup\Gamma_n}(\yb) \quad
\mbox{and}
\quad
\zeta_{\mcF}(\xb)= \displaystyle\sum_{n=1}^{N}\mcX_{\Omega_n\cup\widehat\Gamma_n}(\xb).
\end{equation}
In practice, the piecewise integer-valued function $\zeta_\mcA$ counts how many times, i.e., for how many values of $n$, $\xb$ and $\yb$ belong to $\Omega_n\cup\widehat\Gamma_n\cup\Gamma_n$. The same considerations hold for $\zeta_{\mcF}$. We use these functions to define the bilinear forms and linear functionals associated with each sub-problem; their presence will guarantee that the contributions of overlapping parts of the domains are properly balanced. 

Let $u_{n}$ and $v_{n}$, $n=1,\ldots,N$, be defined on $\Omega_n \cup \widehat{\Gamma}_n\cup\Gamma_n$, we define the subdomain bilinear form and linear functional as follows
\begin{equation}\label{eq:bic} 
\begin{aligned}
&\mathcal{A}_n(u_n,v_n)
\\&
=\int\limits_{\Omega_n\cup\widehat\Gamma_n\cup\Gamma_n}^{}\,\,    
 \int\limits_{\Omega_n\cup\widehat\Gamma_n\cup\Gamma_n}^{}
\hspace{-.2in} 
 \zeta_{\mcA}(\xb,\yb)^{-1} \big(u_n(\yb) - u_n(\xb)\big)\big(v_n(\yb) - v_n(\xb)\big)\gamma(\xb,\yb)d\yb  d\xb, 
\end{aligned}
\end{equation}
where the integrals over $\Gamma_n$ vanish for floating domains, and 
\begin{equation}\label{eq:lfc} 
\mcF_n(v_n) = \int_{\Omega_n\cup\widehat\Gamma_n} \zeta_{\mcF}(\xb)^{-1}v_n(\xb)f(\xb)d\xb.
\end{equation}
We then formulate the {\it domain-decomposition} or multi-domain system as follows: given $f\in W'$, $g\in W_\Gamma$, and a kernel $\gamma$, for $n=1,\ldots, N$, find $u_n\in W_n$ such that
\begin{equation}\label{eq:multi-domain}
\begin{aligned}
&\mathcal{A}_n(u_n,v_n) = \mcF_n(v_n) 
\quad\;\;\;\forall\, v_n\in \widetilde W_n \quad
\mbox{subject to}\\[1mm]
&\left\{\begin{aligned}
 \left(u_n - u_{n'}\right)|_{\widehat\Gamma_{n} \cap \widehat{\Gamma}_{n'}} &= 0 && {\rm for}\,\, n<n'\\
 \left(u_n - g\right)|_{\Gamma_{n}} &= 0,
\end{aligned}\right.
\end{aligned}
\end{equation}
where the last condition only applies to non-floating subdomains. For $n=1,\ldots,N$, the function spaces in \eqref{eq:multi-domain} are defined as follows:
\begin{align}\label{eq:Wn-spaces}
  W_n &= \{w\in L^2 (\Omega_n\cup\widehat\Gamma_n\cup\Gamma_n) : |||w|||_n < \infty\}, \\
  {W_n^0} &= \{w\in W_n : w|_{\Gamma_{n}}=0 \},
\end{align}
where $|||w|||_n^2 =\mathcal{A}_n(w,w) + \|w\|^2_{L^2(\Omega_n\cup\Gamma_n\cup\widehat\Gamma_n)}$, and
\begin{equation}\label{wwtt}
\widetilde W_n = \left\{
\begin{aligned}
&W_n^0 \quad\mbox{if $\Gamma_n\ne\emptyset$}, \,\,\,\,\mbox{(non-floating subdomains)}
\\
&W_n \quad\mbox{if $\Gamma_n=\emptyset$},\,\,\,\,\mbox{(floating subdomains).}
\end{aligned}
\right.
\end{equation}

Paper \cite{Capodaglio2021} proves that the multi-domain system \eqref{eq:multi-domain} is equivalent to the single-domain system \eqref{eq:weak_nonloc}. In the next section we rewrite each sub-problem in \eqref{eq:multi-domain} as an energy-minimization problem and propose a Lagrange-multiplier approach for its solution.

\section{A meshfree substructuring approach}\label{sec:FETI}
Borrowing the idea of finite element tearing and interconnect (FETI) first proposed by \cite{farhat1991method}, we solve the nonlocal domain decomposition problem \eqref{eq:multi-domain} using a Lagrange-multiplier iterative substructuring method. 
We first rewrite problems \eqref{eq:weak_nonloc} and \eqref{eq:multi-domain} as an energy-minimization problem, then introduce a meshfree discretization, and finally propose a solution method.

The energy functional corresponding to \eqref{eq:weak_nonloc} is given by
\begin{align*}
  \mc E(u) &= \frac{1}{2} \mcA(u,u)-\mcF(u) \\
  &=\frac{1}{2}\int_{\Omega\cup\Gamma}  \int_{\Omega\cup\Gamma}\big(u(\yb) - u(\xb)\big)^{2} \gamma(\xb,\yb)d\yb  d\xb - \int_\Omega u(\xb)f(\xb)d\xb
\end{align*}
subject to the Dirichlet volume constraint \(u|_{\Gamma} = g\).

Let a meshfree discretization be given by nodes \(\Lambda=\{\xb_{i}\}_{i=1}^{s}\subset \Omega\) and nodes \(\Lambda^{\Gamma}=\{x_{i}^{\Gamma}\}_{i=1}^{b}\subset \Gamma\).
Then the energy functional can be approximated as
\begin{align*}
  \mcE(u) &\approx \widetilde\mcE(u)\\
          &= \frac{1}{2}\sum_{i=1}^{s}\sum_{j=1}^{s} (  u(\xb_j) - u(\xb_i) )^2 \gamma(\xb_i,\xb_j) \omega_{i}\omega_j \\
          &\quad  + \sum_{i=1}^{s}\sum_{j=1}^{b} (  u(\xb_j^{\Gamma}) - u(\xb_i) )^2 \gamma(\xb_i,\xb_j^{\Gamma}) \omega_{i}\omega_j^{\Gamma} \\
          &\quad + \frac{1}{2}\sum_{i=1}^{b}\sum_{j=1}^{b} (  u(\xb_j^{\Gamma}) - u(\xb_i^{\Gamma}) )^2 \gamma(\xb_i^{\Gamma},\xb_j^{\Gamma}) \omega_{i}^{\Gamma}\omega_j^{\Gamma}  \\
          &\quad - \sum_{i=1}^{s} u(\xb_i)f(\xb_i)\omega_i,
\end{align*}
where \(\{\omega_{i}\}_{i=1}^{s}\) and \(\{\omega_{i}^{\Gamma}\}_{i=1}^{b}\) are numerical quadrature weights.
By taking into account the Dirichlet volume condition on \(\Gamma\), this can equivalently be written as
\begin{equation}
  \widetilde\mcE(\mbf u) = \frac{1}{2}\mbf u^T A \mbf u - \mbf f \cdot \mbf u + const, \\
\label{eqn:global_discrete_matrix_eqns}
\end{equation}
with stiffness matrix \(A\), load vector \(\mbf f\) and solution vector \(\mbf u\), i.e.
\begin{align*}
  [\mbf u]_{i} &= u(\xb_i), & i&\in\Lambda,\\
  [A]_{ij} &=  \delta_{ij}2\sum_{k=1}^{s}\gamma(\xb_{i},\xb_{k})\omega_{i}\omega_{k}-2\gamma(\xb_i,\xb_j) \omega_{i}\omega_j + \delta_{ij}2\sum_{k=1}^{b}\gamma(\xb_{i},\xb_{k}^{\Gamma})\omega_{i}\omega_{k}^{\Gamma}, & i,j&\in\Lambda, \\
  [\mbf f]_{i} &= \omega_{i}f(\xb_{i}) - 2\sum_{j=1}^{b}g(\xb_{j}^{\Gamma})\gamma(\xb_{i},\xb_{j}^{\Gamma})\omega_{i}\omega_{j}^{\Gamma}, &i&\in\Lambda.
\end{align*}
Equivalently, this corresponds to solving the linear system
\begin{equation}
  A \mbf u = \mbf f. \label{eqn:global_eqns_primal}
\end{equation}

In a similar fashion, note that solving problem \eqref{eq:multi-domain} is equivalent to minimizing the energy functionals
\begin{align}
    \mc E_n(u_n) =& \frac{1}{2}\int_{\Omega_n \cup \widehat{\Gamma}_n\cup\Gamma_n} \int_{\Omega_n \cup \widehat{\Gamma}_n\cup\Gamma_n} \zeta_{\mc A}(\xb,\yb)^{-1} (  u_n(\yb) - u_n(\xb) )^2 \gamma(\xb,\yb) d\yb d\xb \notag \\
    & - \int_{\Omega_n\cup\widehat{\Gamma}_n} \zeta_F(\xb)^{-1} u_n(\xb)f(\xb) d\xb   \qquad \qquad \qquad \textrm{for } n = 1, \dots, N,
\end{align}
subject to the nonlocal continuity constraints
\begin{align*}
\left(u_n - u_{n'}\right)|_{\widehat{\Gamma}_n\cap\widehat{\Gamma}_{n'}}&=0 &  \textrm{for } 1\leq n<n'\leq N,
\end{align*}
and Dirichlet volume constraints
\begin{align*}
    \left(u_n - g\right)|_{\Gamma_n}=0.
\end{align*}
Denote in turn $\Lambda_n := \{\xb_{i}^{n}\}_{i=1}^{s_{n}}:= \Lambda\cap(\Omega_{n}\cup\widehat{\Gamma}_{n})$ and $\Lambda_{\Gamma_{n}} := \{\xb_{i}^{\Gamma_{n}}\}_{i=1}^{b_{n}}:= \Lambda\cap\Gamma_{n}$ the nodes of the meshfree discretization of the subdomain $\Omega_n \cup \widehat{\Gamma}_n$ and \(\Gamma_{n}\) respectively.
Then we can discretize the energy functionals and constraints as
\begin{equation}
\left\{
\begin{aligned}
  \mcE_n(u_n) &\approx \widetilde\mcE_n(u_n)\\
              &= \frac{1}{2}\sum_{i=1}^{s_n}\sum_{j=1}^{s_n} \zeta_{\mc A}(\xb_i^{n},\xb_j^{n})^{-1} (  u_n(\xb_j^{n}) - u_n(\xb_i^{n}) )^2 \gamma(\xb_i^{n},\xb_j^{n}) \omega_{i}^{n}\omega_j^{n} \\
              &\quad + \sum_{i=1}^{s_n}\sum_{j=1}^{b_n} \zeta_{\mc A}(\xb_i^{n},\xb_j^{\Gamma_{n}})^{-1} (  u_n(\xb_j^{\Gamma_{n}}) - u_n(\xb_i^{n}) )^2 \gamma(\xb_i^{n},\xb_j^{\Gamma_{n}}) \omega_{i}^{n}\omega_j^{\Gamma_{n}} \\
              &\quad + \frac{1}{2}\sum_{i=1}^{b_n}\sum_{j=1}^{b_n} \zeta_{\mc A}(\xb_i^{\Gamma_{n}},\xb_j^{\Gamma_{n}})^{-1} (  u_n(\xb_j^{\Gamma_{n}}) - u_n(\xb_i^{\Gamma_{n}}) )^2 \gamma(\xb_i^{\Gamma_{n}},\xb_j^{\Gamma_{n}}) \omega_{i}^{\Gamma_{n}}\omega_j^{\Gamma_{n}} \\
              &\quad - \sum_{i=1}^{s_n}  \zeta_F(\xb_i^{n})^{-1}u_n(\xb_i^{n})f(\xb_i^{n})\omega_i^{n}, \\
  u_n(\mbf x) &= u_{n'}(\mbf x) \qquad \mbf x \in (\widehat{\Gamma}_n\cap\widehat{\Gamma}_{n'}) \cap (\Lambda_n \cap \Lambda_{n'}),
\end{aligned} \right.
\label{eq:discrete_eqns}
\end{equation}
where $\{\omega_{i}^{n}\}_{i=1}^{s_n}$ and $\{\omega_{i}^{\Gamma_{n}}\}_{i=1}^{s_n}$ are a set of numerical quadrature weights.
Taking into account the Dirichlet volume condition on \(\Gamma_n\), we simplify the expressions in \eqref{eq:discrete_eqns} with the following matrix notation:
\begin{equation}
\left\{
\begin{aligned}
&  \widetilde\mcE_n(\mbf u_n) = \frac{1}{2} \mbf u_n^T A_n \mbf u_n - \mbf f_n \cdot \mbf u_n  + const, \\
&  M_{nn'} \mbf u_n  - M_{n'n} \mbf u_{n'}=\mbf 0, \quad\textrm{for } n < n',
\end{aligned} \right.
\label{eqn:discrete_matrix_eqns}
\end{equation}
where $\mbf u_n \in \mathbb{R}^{s_n}$ is the vector with $[\mbf u_n]_{i}=u_n(\xb_i^{n})$ and the entries of stiffness matrix $A_n\in \mathbb{R}^{s_n\times s_n}$ and of the vector $\mbf f_n \in \mathbb{R}^{s_n}$ are given by
\begin{align*}
  [A_n]_{ij} &=  
  \delta_{ij}2\sum_{k=1}^{s}\zeta_{\mc A}(\xb_i^{n},\xb_k^{n})^{-1}\gamma(\xb_{i}^{n},\xb_{k}^{n})\omega_{i}^{n}\omega_{k}^{n} \\
  & \quad - 2\zeta_{\mc A}(\xb_i^{n},\xb_j^{n})^{-1}\gamma(\xb_i^{n},\xb_j^{n}) \omega_{i}^{n}\omega_j^{n}  \\ 
  & \quad + \delta_{ij}2\sum_{k=1}^{b}\zeta_{\mc A}(\xb_i^{n},\xb_k^{\Gamma_n})^{-1}\gamma(\xb_{i}^{n},\xb_{k}^{\Gamma_n})\omega_{i}^{n}\omega_{k}^{\Gamma_n}, & i,j&\in\Lambda_{n}, \\
  [\mbf f_n]_{i}  &= \zeta_F(\xb_i^{n})^{-1}f(\xb_i^{n})\omega_i^{n} - \sum_{j=1}^{b}\zeta_{\mc A}(\xb_i^{n},\xb_j^{\Gamma_n})^{-1}g(\xb_{j}^{\Gamma_n})\gamma(\xb_{i}^{n},\xb_{j}^{\Gamma})\omega_{i}^{n}\omega_{j}^{\Gamma_n}, &i&\in\Lambda_{n}.
\end{align*}
Let $p_{nn'}$ denote the number of degrees of freedom in $\Lambda_n\cap \Lambda_{n'}$.
Then, the matrix $M_{nn'} \in \mathbb{R}^{p_{nn'}\times s_n}$ selects the nodes of subdomain $n$ that are in the overlapping region with subdomain $n'$.

To determine the solution to the constrained minimization problem \eqref{eqn:discrete_matrix_eqns} with a Lagrange-multiplier approach, we introduce the multiplier $\bs \lambda\in\mbR^{q}$, where $q$ is number of the nonlocal continuity constraints in total, $q := \sum_{n=1}^{N} \sum_{n'=n+1}^{N} p_{nn'}$. Following the procedure and the notation illustrated in \cite{mathew2008}, Chapter 4, for Lagrange-multiplier based substructuring algorithms, we rewrite the multi-domain system \eqref{eqn:discrete_matrix_eqns} as the following system of global extended matrix equations:
\begin{align}
  \left(
  \begin{array}{cc}
    A_{\epsilon\epsilon}&M^{T}\\
    M&0
  \end{array}
       \right)
       \left(
       \begin{array}{c}
         \mbf u_{\epsilon}\\
         \bs \lambda
       \end{array}
\right) = \left(
       \begin{array}{c}
         \mbf f_{\epsilon}\\
         \mbf 0
       \end{array}
\right)
\label{eqn:global_eqns}
\end{align}
where the {\it extended} stiffness matrix $A_{\epsilon\epsilon}$, displacement $\ub_\epsilon$ and right-hand side $\mbf f_\epsilon$ are defined as
\begin{align*}
A_{\epsilon\epsilon} &= \left( \begin{array}{ccccc}
A_1 &  &  &  &  \\
 & A_2 &  &  &  \\
 &  & \ddots & &  \\
 &  &  &  & A_{N} \\
\end{array} \right), 
&\ub_{\epsilon} &= \left( \begin{array}{c}
\ub_1   \\
\ub_2  \\
\vdots   \\
\ub_{N} \\
\end{array} \right), & \mbf f_{\epsilon} &= \left( \begin{array}{c}
\mbf f_1   \\
\mbf f_2  \\
\vdots   \\
\mbf f_{N} \\
\end{array} \right).
\end{align*}
Equivalently, \(A_{\epsilon\epsilon}\) can be written as
\begin{equation}
  A_{\epsilon\epsilon} = \sum_{n=1}^{N} (R^P_{n})^TA_{n}R^P_{n}, \label{eqn:distributedA}
\end{equation}
where $R^P_{n} = \left(\bs 0  \ldots I_{s_{n}} \ldots \bs 0\right)$ is the restriction from the global extended primal space, and \(I_{k}\) is the identity matrix of size \(k\).

Let $s_\epsilon := \sum_{n=1}^{N} s_n$ denote the number of degrees of freedom of the extended system $A_{\epsilon\epsilon}$.
Then the matrix $M\in\mbR^{q\times s_\epsilon}$ is defined in such a way that the continuity constraints in the second equation in \eqref{eqn:discrete_matrix_eqns} are satisfied.
Since the stiffness matrices of the floating subdomains are singular, the resulting extended stiffness matrix $A_{\epsilon\epsilon}$ is singular if the set of floating subdomains is non-empty. Let the matrix $Z_n \in \mathbb{R}^{s_n \times z_n}$ denote the matrix whose columns span the nullspace of the subdomain stiffness matrix $A_n$.
If subdomain $n$ is a non-floating subdomain, let $Z_n$ be a column vector of all zeros and $z_n=1$.
Then we define a block matrix $Z \in \mathbb{R}^{s_\epsilon \times z}$ with $z = \sum_{n=1}^{N} z_n$ as
\begin{align}\label{eqn:distributedZ}
  Z&= \sum_{n=1}^{N} (R^P_{n})^T Z_{n}R^C_{n}
\end{align}
with
$R^C_{n} = \left(\bs 0  \ldots I_{z_{n}} \ldots \bs 0 \right)^{T}$ the restriction from the global coarse index space. The matrix $M$ is given by
%
\begin{equation}
    M=\sum_{n=1}^{N} (R^D_n)^T M_n R^P_n,
    \label{eqn:distributedM}
\end{equation}
where
\begin{align*}
 (R^D_n)^T &= \left( \begin{array}{cccccccc}
 \bs 0 &  &   &   & & & & \\
 I_{p_{n1}} & &  &  &  & & & \\
  \bs 0 &  \bs 0 & &  & & & & \\
&  I_{p_{n2}} & & & & & & \\
&  \bs 0 & \ddots & & & & & \\
&   &  & \bs 0 & & & & \\
&   &  & I_{p_{n,n-1}} & & & &\\
&   &  & \bs 0 & & & &\\
&   &  &  & I_{p_{n,n+1}} & & &\\
&   &  &  &  &  \ddots & &\\
&   &  &  &  &   & I_{p_{n,N}} &\\
&   &  &  &  &   & & \bs 0\\
\end{array} \right), \; \notag \\
    M_n &= \left( \begin{array}{c}
      -M_{n1}\\
      \vdots\\
      -M_{n,n-1}\\
      M_{n,n+1}\\
      \dots\\
      M_{n,N}\\
    \end{array} \right).
\end{align*}
Hence, $R^D_n \in \mathbb{R}^{q_n \times q}$ is the restriction from the global dual index space and $M_n\in\mathbb{R}^{q_n\times s_n}$, where $q_n=\sum_{n'=1}^{N}p_{n,n'}$.

With implementation in mind, we note that matrices \(A_{n}\), \(M_{n}\) and \(Z_{n}\), as well as vectors \(\mbf u_{n}\) and \(\mbf f_{n}\)  are local to subdomain \(n\).
Matrices \(R^P_{n}\), \(R^C_{n}\) and \(R^D_{n}\) describe the injection into global index spaces and hence encode information exchange between subdomains.

\subsection{Reduced system for the Lagrange multiplier and its solution}
In this section we show how to derive a reduced system for the Lagrange multiplier $\bs\lambda$ and describe an iterative method for its numerical solution. 

By projecting the first equation of \eqref{eqn:global_eqns} onto the nullspace, the following compatibility condition is derived:
\begin{align}
  Z^T(\mbf f_\epsilon - M^T \bs\lambda ) = \mbf 0
  \quad \Leftrightarrow \quad G^T \bs\lambda = \mbf g,
\end{align}
where $Z$ is the matrix that spans the nullspace of $A_{\epsilon\epsilon}$, $G := MZ$ and $\mbf g := Z^T \mbf f_\epsilon$. By using the condition $G^T {\bs\lambda} = \mbf g$, the generalized solution to the first block row in \eqref{eqn:global_eqns} can be written as
\begin{align}
    \mbf u_\epsilon = A^\dagger_{\epsilon\epsilon}(\mbf f_\epsilon - M^T \bs\lambda) + Z \bs\alpha,
    \label{eqn:general_sol_uepsilon}
\end{align}
where $\bs\alpha \in \mathbb{R}^d$ is yet to be determined vector and the symbol $\dagger$ denotes the \emph{Moore-Penrose pseudoinverse}. By substituting \eref{general_sol_uepsilon} into the continuity constraints equation in \eref{global_eqns} we obtain
\begin{align}
    K \bs\lambda - G \bs\alpha = \mbf e,
\end{align}
where $K = M A_{\epsilon\epsilon}^\dagger M^T $ and $\mbf e = MA_{\epsilon\epsilon}^\dagger \mbf f_\epsilon$. The term $G\bs\alpha$ can be further eliminated by applying the projection $P_0 := I - G(G^TG)^\dagger G^T$. Therefore, together with the compatibility equation, the reduced system for the Lagrange multiplier $\bs\lambda$ can be written as
\begin{align}
  \left\{
  \begin{array}{rcl}
    P_{0}K\bs \lambda&=&P_{0} \mbf e\\
    G^{T}\bs \lambda&=&\mbf g.
  \end{array}\right.
  \label{eqn:reduced_system}
\end{align}
To solve the reduced system $\eref{reduced_system}$, we utilize the preconditioned iterative parallel algorithm reported in Algorithm \ref{algo:parallelAlgorithm}. We describe the choice of the preconditioner in Section \ref{sec:prec} and the extension of the standard iterative solver to the parallel setting in Section \ref{sec:distributed}.

\subsubsection{A preconditioner based on the Schur complement}\label{sec:prec}
We exploit the diagonal block structure of $A_{\epsilon\epsilon}$ and equation \eref{distributedM} to write the matrix $K$ as
\begin{align}
  K &= M A_{\epsilon\epsilon}^{\dagger}M^{T} \notag \\
    &= \left(\sum_{n=1}^{N} (R^D_n)^T M_n R^P_n\right) \left(\sum_{n=1}^{N} (R^P_n)^T A_{n}^{\dagger}R^P_{n}\right) \left(\sum_{n=1}^{N} (R^P_n)^T M_n^{T} R^{D}_n\right) \notag \\
    &= \sum_{n=1}^{N} (R^D_n)^T M_n A_n^\dagger M_n^T R^{D}_n, \label{eq:localK}
\end{align}
where we have used \eqref{eqn:distributedA}, \eqref{eqn:distributedM} and the fact that \(R^P_{n}(R_{n'}^{P})^T=\delta_{nn'}I_{s_{n}}\).
Without loss of generality, we assume that the columns of each matrix $M_n$ are ordered in such a way that the interior nodes precede the interface nodes on \(\widehat{\Gamma}_{n}\), i.e., 
\begin{align*}
    M_n = [M_n^{(\Omega_n)} \quad M_n^{(\widehat{\Gamma}_n)}] = [0 \quad M_n^{(\widehat{\Gamma}_n)}],
\end{align*}
where $M_n^{(\Omega_n)} = 0$ because the continuity constraints only involve boundary nodes. Similarly, matrix $A_n$ can also be written as 
\begin{align*}
A_n &= \left( \begin{array}{cc}
A_n^{(\Omega_n\Omega_n)} A_n^{(\Omega_n \widehat{\Gamma}_n)} \\
A_n^{(\widehat{\Gamma}_n \Omega_n)} A_n^{(\widehat{\Gamma}_n \widehat{\Gamma}_n)}
\end{array} \right).
\end{align*}
Using \eqref{eq:localK} the matrix $K$ can then be written as
\begin{align*}
    K &= \sum_{n=1}^{N} (R^D_n)^T[0 \quad M_n^{(\widehat{\Gamma}_n)}] A_n^\dagger [0 \quad M_n^{(\widehat{\Gamma}_n)}]^T R^{D}_n \notag \\
    &= \sum_{n=1}^{N} (R^D_n)^T M^{(\widehat{\Gamma}_n)}_n S_n^\dagger M^{(\widehat{\Gamma}_n)^T}_n R^{D}_n,
\end{align*}
where $S_n$, the \emph{Schur complement} with respect to boundary nodes of $A_n$, is defined as
\begin{align*}
    S_n = A_n^{(\widehat{\Gamma}_n \widehat{\Gamma}_n)} - A_n^{(\widehat{\Gamma}_n \Omega_n)} A_n^{(\Omega_n \Omega_n)^{-1}} A_n^{(\Omega_n \widehat{\Gamma}_n)}.
\end{align*}
Since $(M_n^{(\widehat{\Gamma}_n)})^T R^D_n$ restricts, up to a sign change, to the Lagrange multipliers associated with subdomain $n$, the structure of \(K\) motivates the following Schur complement-based preconditioner $Q$ \citep{mathew2008}: 
\begin{align}
    Q = \sum_{n=1}^{N} (R^D_n)^T M^{(\widehat{\Gamma}_n)}_n S_n M^{(\widehat{\Gamma}_n)^T}_n R^{D}_n.
\end{align}
This preconditioner is commonly named the \emph{Dirichlet preconditioner}.

\subsubsection{The \textit{distributed} projected gradient algorithm}\label{sec:distributed}

For an arbitrary global vector $\bs \lambda$, the matrix-vector multiplication \(\bs \mu=K\bs \lambda\) with the global matrix $K$ can be written using \eqref{eq:localK} as
\begin{align*}
  \bs \mu =& \sum_{n=1}^{N} (R^D_n)^T M_n A_n^\dagger M_n^T R^{D}_n \bs \lambda,
\end{align*}
or equivalently
\begin{align*}
    R^{D}_{n'} \bs \mu =& \sum_{n=1}^{N} (R^{D}_{n'} (R^D_n)^T) M_n A_n^\dagger M_n^T (R^{D}_n \bs \lambda).
\end{align*}
Moreover, if we set \(\bs \lambda_{n}:=R^{D}_n \bs \lambda\) and \(\bs \mu_{n}:=R^{D}_n \bs \mu\) for the local parts of \(\bs \lambda\) and \(\bs \mu\) respectively, and $K_n := M_n A^\dagger_n M^T_n$, we obtain
\begin{align}
    \bs \mu_{n'} =& \sum_{n=1}^{N}(R^{D}_{n'} (R^D_n)^T) K_n \bs \lambda_{n}.
    \label{eqn:dist_K}
\end{align}
Except for neighboring subdomains \(n\) and \(n'\), the product \(R^{D}_{n'} (R^D_{n})^T\) is zero.
Hence, the sum can be understood in terms of a halo exchange.

Therefore, the global matrix-vector product can be expressed in terms of local matrix-vector product with \( K_n \) and a round of neighborhood communication.
Similarly, the action of the preconditioner \(Q\) reduces to a multiplication with \( Q_n := M_n^{(B)} S_n M^{(B)^T}_n\) and accumulation over neighboring subdomains.

Due to \eqref{eqn:distributedM}, \eqref{eqn:distributedZ} and \(R^P_{n} (R^{P}_{n'})^T = \delta_{nn'}I_{s_{n}}\), we have that
\begin{align*}
  G &= MZ = \sum_{n=1}^{N} (R^D_{n})^T M_{n}Z_{n} R^C_{n},
\end{align*}
and hence
\begin{align*}
  G^{T}\bs \lambda &= \sum_{n=1}^{N} (R^{C}_{n})^T Z_{n}^{T} M_{n}^{T}(R^{D}_{n} \bs \lambda) \\
                   &= \sum_{n=1}^{N} (R^{C}_{n})^T Z_{n}^{T} M_{n}^{T} \bs \lambda_{n}.
\end{align*}
We elect to solve the coarse-grid problem \(G^{T}G\) redundantly on all ranks. 
Alternatively, one could solve it only on a subset of the subdomains, or even a single rank, and then scatter the solution back to all subdomains.
Therefore, the sum in the above expression can be interpreted as a global reduction of the results of the local matrix-vector products \(Z_{n}^{T}M_{n}\bs\lambda_{n}\).

Thus, the action of the projection \(P_{0}\) on a vector \(\bs \lambda\) can be written as
\begin{align*}
  \bs \mu &:= P_{0}\bs \lambda \\
          &= \bs \lambda - G (G^{T}G)^{\dagger} G^{T}\bs\lambda \\
          &= \bs\lambda - \left(\sum_{n=1}^{N} (R^{D}_{n})^T M_{n}Z_{n} R^C_{n}\right) (G^{T}G)^{\dagger} \left(\sum_{n=1}^{N} (R^{C}_{n})^T Z_{n}^{T}M_{n}^{T} \bs \lambda_{n}\right)
\end{align*}
and, consequently,
\begin{align}
  \bs \mu_{n'} &= \bs\lambda_{n'} - \left(\sum_{n=1}^{N} R^{D}_{n'} (R^{D}_{n})^T M_{n}Z_{n} R^C_{n}\right) (G^{T}G)^{\dagger} \left(\sum_{n=1}^{N} (R^{C}_{n})^T Z_{n}^{T}M_{n}^{T} \bs \lambda_{n}\right).
  \label{eqn:dist_P0}
\end{align}
In order to avoid another halo exchange during the repeated apply of \(P_{0}\), we precompute the action of
\begin{align*}
  \widetilde{G}_{n'} &:= R^{D}_{n'} G = \sum_{n=1}^{N} R^{D}_{n'} (R^{D}_{n})^T M_{n}Z_{n} R^C_{n}.
\end{align*}
In our distributed preconditioned projected gradient algorithm we use the notation $\xleftarrow{R^D_n}$ and $\xleftarrow{R^D_n}$ to represent the operations of $\sum_{n=1}^{N}R^{D}_{n'} (R^D_n)^T$ and $\sum_{n=1}^{N} (R^{C}_{n})^T$ as described in \eref{dist_K} and \eref{dist_P0} respectively.
In our implementation these halo exchanges are realized using \texttt{MPI\_Irecv} and \texttt{MPI\_Isend}.
Furthermore, inner products of distributed vectors need to be computed. 
Using \texttt{MPI\_Allreduce} in practice, the global reductions are denoted by $\xleftarrow{\textrm{reduction}}$ in Algorithm~\ref{algo:parallelAlgorithm}.

\begin{algorithm}
    
    \For{each processor $n$}{
    \emph{assemble local matrices $A_n,\mbf f_n,M_n,Z_n$ and preconditioner $Q_n$}\;
    \emph{assemble and compute $(G^T G)^\dagger$ and $\widetilde{G_n}$}\;
    $g_{n} = Z_n^T \mbf f_n, \quad \mbf e_n = M_n A^\dagger_n f_n, \quad G_n = M_n Z_n, \quad K_n = M_n A_n^\dagger M_n^T $\;
    \emph{Calculate distributed initial guess }$\bs \lambda_n = \widetilde{G}_n(G^T G)^\dagger g_n $\;
    $\mbf x_n = K_n \bs \lambda_n ,\quad \mbf x_n \xleftarrow{R^D_n}  \mbf x_{n}$ \;
    $\mbf r_n = \mbf e_n - \mbf x_n, \quad \mbf z_n = Q \mbf r_n, \quad \mbf z_n \xleftarrow{R^D_n}  \mbf z_{n}$ \;
    $\mbf z'_n = G_n^T \mbf z_n, \quad \mbf z'_n \xleftarrow{R^C_n} \mbf z'_n, \quad \mbf y_n = \mbf z_n - \widetilde{G}_n(G^TG)^\dagger \mbf z'_n$ \;
    $\mbf p_n = \mbf y_n, \quad \alpha \xleftarrow{\textrm{reduction}} \mbf r^T_n \mbf y_n;$
    }
    \Repeat{$r^T y$ is small enough}{
    \For{each processor $n$}{
    $\mbf x_n = K_n \mbf p_n,\quad \mbf x_n \xleftarrow{R^D_n} \mbf x_{n}$\;
    $\mbf x'_n = G_n^T \mbf x_n, \quad \mbf x'_n \xleftarrow{R^C_n} \mbf x'_n, \quad \mbf x_n = \mbf x_n - \widetilde{G}_n(G^TG)^\dagger \mbf x'_n$ \;
    $\beta \xleftarrow{\textrm{reduction}} \mbf p^T_n \mbf x_n, \quad \mbf r^+_n = \mbf r_n - \frac{\alpha}{\beta} \mbf p_n, \quad \bs \lambda_n = \bs\lambda_n + \frac{\alpha}{\beta} \mbf p_n $\ ;
    
    
    $\mbf z_n = Q_n \mbf r^+_n, \quad \mbf z_n \xleftarrow{R^D_n} \mbf z_{n}$ \;
    
    $\mbf z'_n = G_n^T \mbf z_n, \quad \mbf z'_n \xleftarrow{R^C_n} \mbf z'_n, \quad \mbf y^+_n = \mbf z_n - \widetilde{G}_n(G^TG)^\dagger \mbf z'_n$ \;

    
    $\alpha^+ \xleftarrow{\textrm{reduction}} (\mbf r^+_n)^T \mbf y^+_n ,  \quad \mbf p_n = \mbf y_n^+ + \frac{\alpha^+}{\alpha}\mbf p_n$\
    
    $\mbf y_n \leftarrow \mbf y_n^+, \quad \mbf r_n\leftarrow \mbf r_n^+, \quad \alpha\leftarrow \alpha^+ $
    }
    }
    \caption{Distributed Preconditioned Projected CG}
    \label{algo:parallelAlgorithm}
\end{algorithm}

\section{Numerical experiments}
\label{sec:numer-exper}

In this section we illustrate the efficiency and scalability of our approach with several two-dimensional numerical tests. We first provide details regarding the problem setting and discretization, the hardware used in our simulations, and the expected scalability properties of our algorithm.

\paragraph{Problem setting}
Let $\Omega=[0,1]^2$ and let $B_\delta(\xb)=\{|\xb-\yb|_{\ell^\infty}\leq\delta\}$, i.e. we consider square nonlocal neighborhoods\footnote{Note that even though Euclidean neighborhoods are the standard choice in nonlocal modeling, some recent works have considered the use of square neighborhoods, for which the associated nonlocal problem is well posed \citep{Capodaglio2020,Xu2020learning}.}. In all our tests we consider the kernel $\gamma(\xb,\yb)=C|\xb-\yb|^{-1}$ with scaling constant $C\sim\delta^{-3}$, and the manufactured solution $u(\xb)=\xb_1^2+\xb_2^2$ for which the corresponding forcing term is given by $f(\xb)=-4$. 

For the meshfree discretization described in Section \ref{sec:FETI}, we consider \(L\times L\), \(L\in\mathbb{N}^+\), uniformly spaced particles in $\Omega$, resulting in a particle spacing of \(h=L^{-1}\).
Each particle interacts with the square neighborhood $B_\delta(\xb)$, of \((2m+1)\times(2m+1)\) particles, where \(m\in\mathbb{N}^+\), i.e. \(\delta=m h=m/L\).
In what follows, we always choose a zero initial guess and solve the problem in primal variables \eqref{eqn:global_eqns_primal} to a relative tolerance of 1e$-5$.
The solve times discussed below exclude the setup costs.

\paragraph{Computational resources}
All experiments were performed on the Solo cluster at Sandia National Laboratories \cite{sandiaHPC}.
Solo consists of 374 nodes with dual socket 2.1 GHz Intel Broadwell CPUs.
Each node has 36 cores and 128 GB of RAM.
Solo uses an Intel Omni-Path interconnect.

\paragraph{Expected scalability properties}
We distribute the problem over \(p\times p\) subdomains, so that each subdomain consists of \(s_{i}=(L/p)^{2}\) particles.
Since every particle interacts with \((2m+1)^{2}\) particles, the subdomain matrices \(A_{i}\) have
\begin{align*}
  \operatorname{nnz}(A_{i}) \sim \left(\frac{Lm}{p}\right)^{2} \sim \left(\frac{L^{2}\delta}{p}\right)^{2}
\end{align*}
entries.
The condition number \(\kappa\) of the global linear system \(A\) defined in \eqref{eqn:global_discrete_matrix_eqns} scales like
\begin{align*}
  \kappa \sim \delta^{-2},
\end{align*}
as can be seen by comparing this particular case of a meshfree discretization with finite elements on a uniform mesh \cite{aksoylu2010results,aksoylu2011variational,aksoylu2014conditioning}.

Therefore, we expect
\begin{align}
  I_{\text{Krylov}}\sim\sqrt{\kappa} \sim 1/\delta\sim \frac{L}{m} \label{eq:itsKrylov}
\end{align}
iterations when using a Krylov method to solve the equations involving \(A\).
The cost per iteration, however, will be proportional to the cost of a local matrix-vector product, i.e. \(\operatorname{nnz}(A_{i})\), so that overall solve time \(T_{\text{Krylov}}\) of an unpreconditioned Krylov method such as the conjugate gradient method will scale like
\begin{align}
  T_{\text{Krylov}} \sim I_{\text{Krylov}} \operatorname{nnz}(A_{i}) &\sim \left(\frac{L^{2}}{p}\right)^{2} \delta \sim \left(\frac{L}{p}\right)^{2} Lm. \label{eq:timeKrylov}
\end{align}

We are contrasting this scaling with the proposed domain decomposition solvers.
Solve times are dominated by the projected CG iteration.
Each matrix-vector product involving the operator \(P_{o}K\) involves the parallel solution of one linear system per subdomain.
We are employing a direct solver in this step, meaning that the cost per operator apply will at most be
\begin{align*}
  s_{i}^{2}\sim \left(\frac{L}{p}\right)^{4}.
\end{align*}
It will transpire from numerical experiments that the DD methods display a number of iterations that admits a bound \(I_{\text{DD}}(L/p)\), where \(I_{\text{DD}}(\cdot)\) is monotonically increasing in its argument.
This leads us to expect overall solve times proportional to
\begin{align}
  T_{\text{DD}} \sim~ I_{\text{DD}}(L/p) \left(\left(\frac{L}{p}\right)^{4} + \log p\right). \label{eq:timeDD}
\end{align}
Here, the additional \(\log p\) contribution is due to the projection step, which involves global communication.

Further implementation details of the methods such as the used solvers will be discussed below.

\subsection{Weak scaling with constant \texorpdfstring{\(m\)}{m}}\label{sec:weak-scaling}

\begin{figure}
  \centering
  \includegraphics{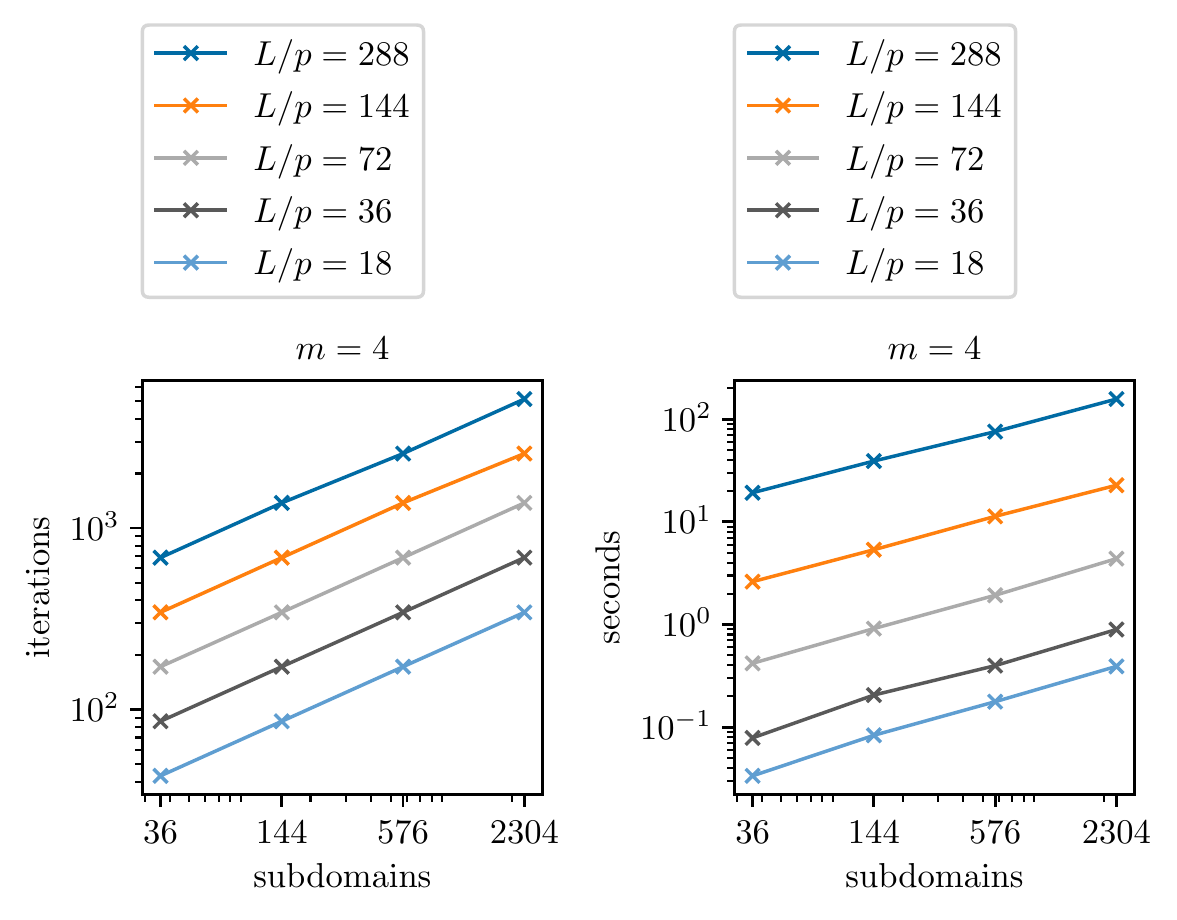}
  \caption{Weak scaling of CG with constant \(m=4\) and \(p\sim L\).}
  \label{fig:weakScalingCG}
\end{figure}

In a first weak scaling experiment, we fix the number of particle interactions \(m=4\) and vary \(L\) and \(p\) proportionally, resulting in constant memory footprint per subdomain.
We first solve the resulting (single-domain) linear equations using the conjugate gradient method, since it is a commonly exercised solver for nonlocal problems, e.g.\ \cite{peridigm}, and therefore a reasonable baseline for comparison.  Here, we distribute the rows of the global matrix according to the non-overlapping partition of \(\Omega\cup\Gamma\).  Timings and iteration counts are displayed in Figure~\ref{fig:weakScalingCG}.
As expected from equations \eqref{eq:itsKrylov} and \eqref{eq:timeKrylov}, it can be observed that the number of iterations and, consequently, the overall solve times scale like \(I_{\text{Krylov}} \sim T_{\text{Krylov}} \sim L \sim p\).
In particular, this shows that a Krylov method without adequate preconditioning will not result in a scalable solver.

We next turn to the proposed substructuring solver.
The subdomain problems are factored using a Cholesky decomposition, provided by CHOLMOD \cite{ChenDavisEtAl2008_Algorithm887}.
\begin{figure}
  \centering
  \includegraphics{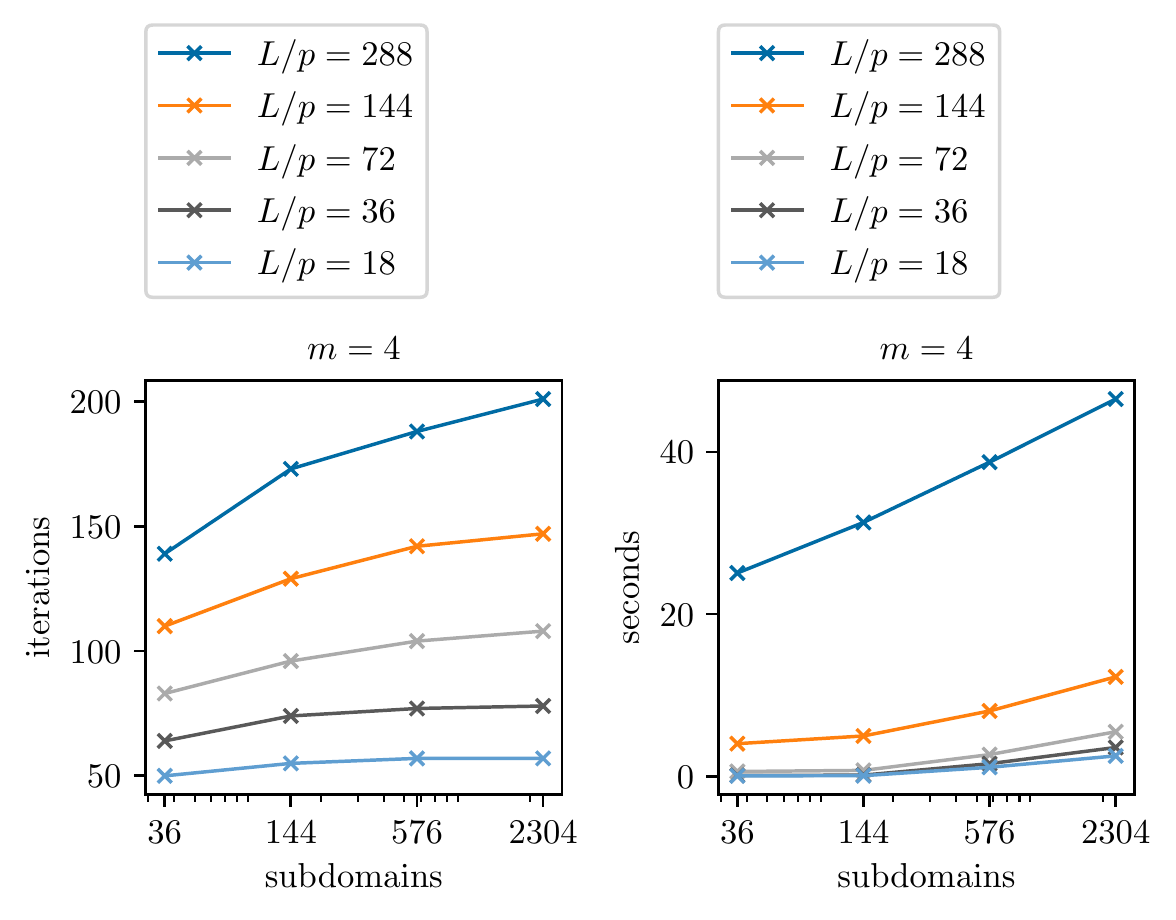}
  \includegraphics{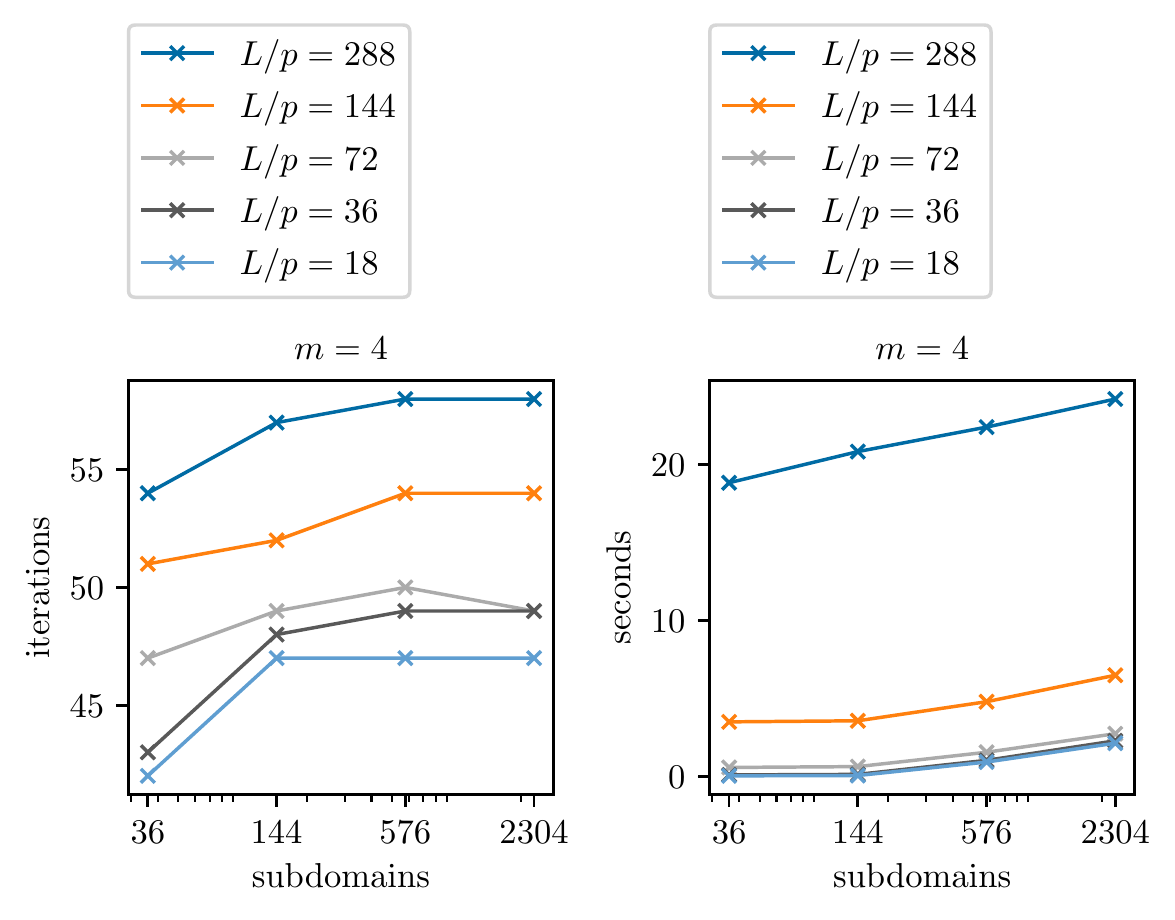}
  \caption{
    Weak scaling of the substructuring solver with constant \(m=4\) and \(p\sim L\).
    \emph{Top:} No preconditioner, Cholesky solver in the operator.
    \emph{Bottom:} Dirichlet preconditioner, Cholesky solvers in the operator and preconditioner.
  }
  \label{fig:weakScalingCholesky}
\end{figure}
We show the number of iterations of the projected CG method and overall solve times in the top row of Figure~\ref{fig:weakScalingCholesky}.
The number of iterations appears to approach an asymptotic value that depends on \(L/p\).
Although the number of subdomains increases by a factor of \(64\), the solve time increases at most by a factor of 1.85, which is in line with the predicted behavior in \eqref{eq:timeDD}.

Next, we precondition the reduced system with the Dirichlet preconditioner described above.
We also use a Cholesky factorization for the involved interior solves.
Again, iterations and solve times are shown in the bottom row of Figure~\ref{fig:weakScalingCholesky}.
The number of iterations varies between 42 and 58, and is asymptotically independent of the problem size and the number of subdomains.
Similarly, the overall solve time is close to constant as the problem size is increased, but depends on the size of the subdomain problems.
Clearly, use of the Dirichlet preconditioner leads to improvements in terms of iteration count and solve time compared to the unpreconditioned DD solver.

\begin{figure}
  \centering
  \includegraphics{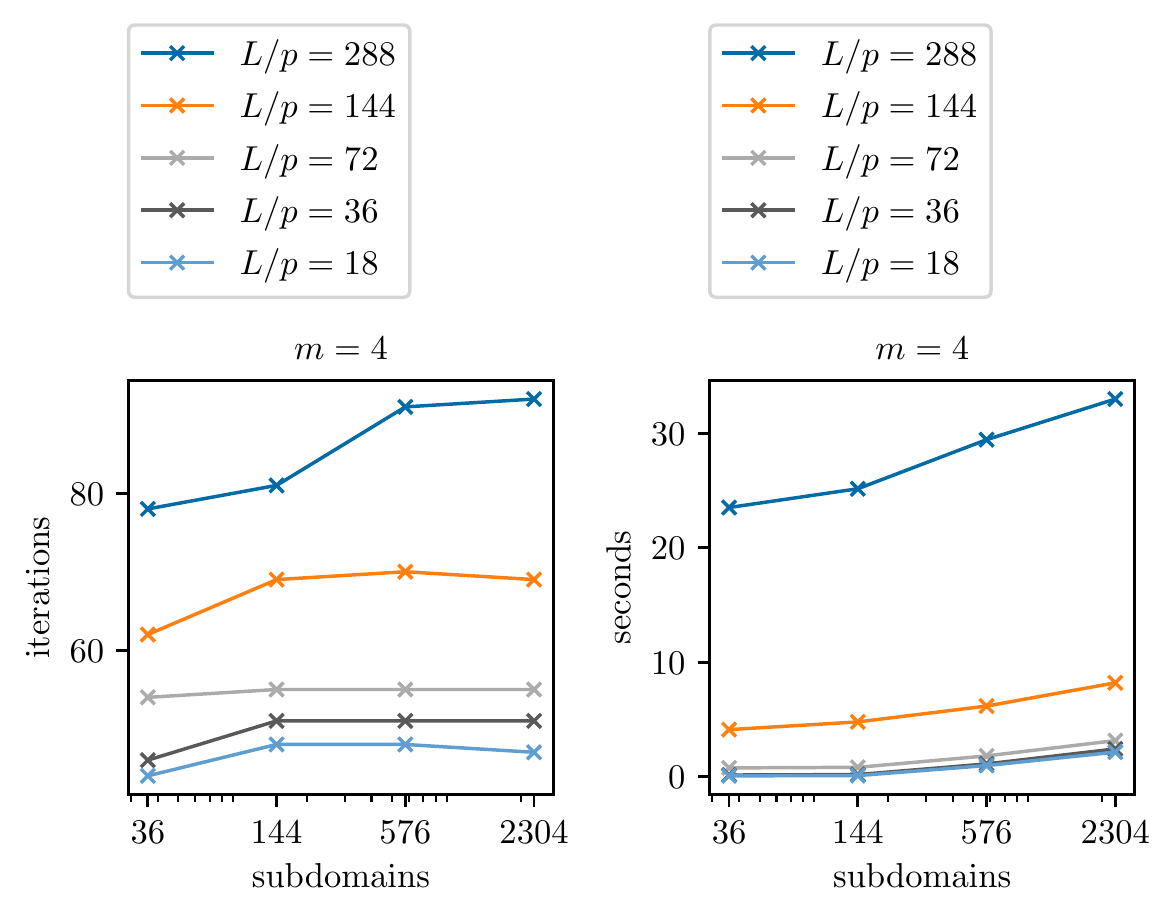}
  \includegraphics{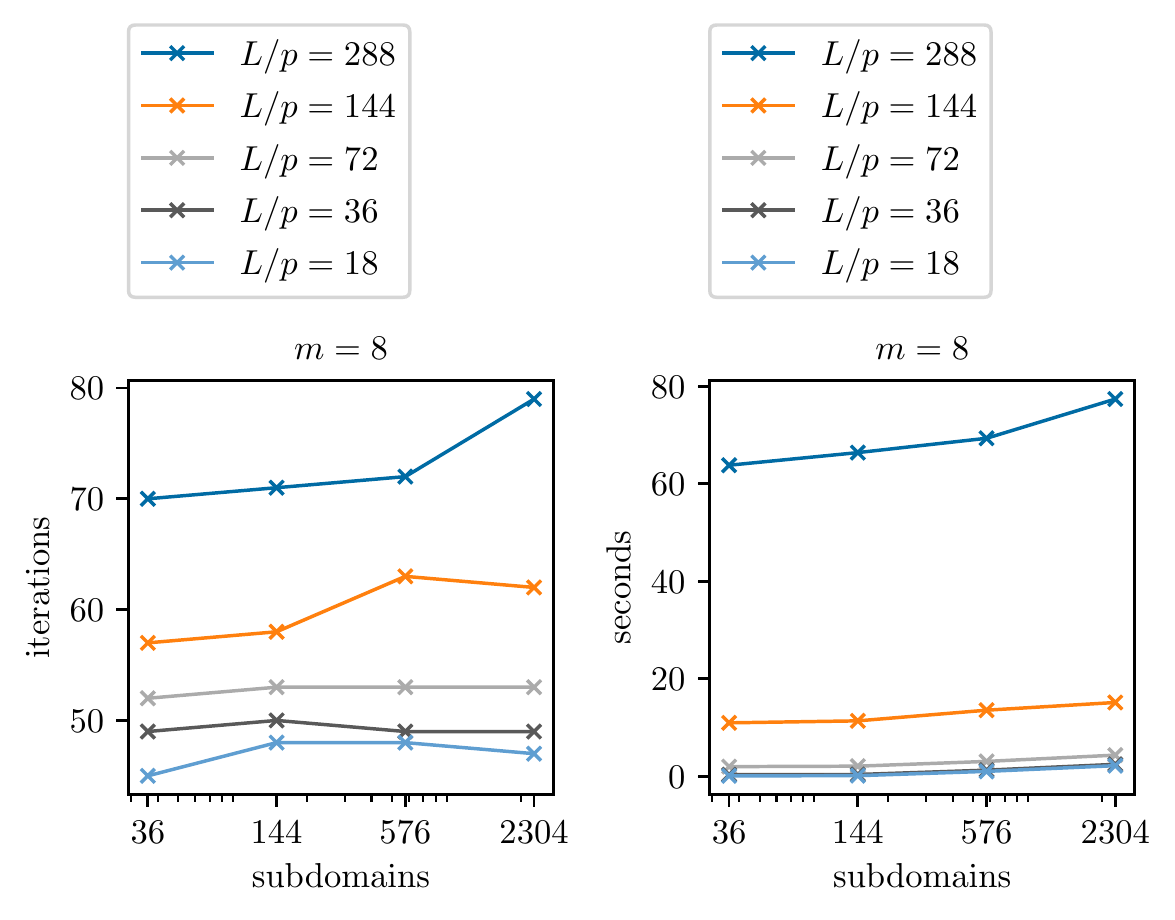}
  \caption{Weak scaling of the substructuring solver with constant \(m\) (\emph{top:} \(m=4\), \emph{bottom:} \(m=8\)) and \(p\sim L\) and Cholesky solver in the operator and 5 iterations of CG in the Dirichlet preconditioner.}
  \label{fig:weakScalingCGDirichletPrec}
\end{figure}

The memory cost of computing two Cholesky factorizations can be significant, in particular for large number of particle interactions \(m\).
It is therefore natural to attempt to reuse the factorization of the interior when solving the subdomain problem, or to replace the direct solvers with an iterative method.
While using an inexact solve for the subdomains within the operator can easily lead to difficulties with regards to convergence towards the correct solution, we can replace the solve on the subdomain interior within the Dirichlet preconditioner with 5 iterations of CG and repeat the weak scaling experiment.
The results in Figure~\ref{fig:weakScalingCGDirichletPrec} show a minor slow-down in the solve, but the overall scaling of the method is unaffected.
We also show timings for \(m=8\).

By selecting the problem with 2304 subdomains and \(L/p=288\) for comparison between the four presented solver options, we observe that the CG solve took 157.16 seconds, while the domain decomposition solver takes 46.5s without preconditioner, 24.2s with Cholesky-Dirichlet preconditioner and 33s with CG-Dirichlet preconditioner.
Based on the scaling estimates, speedups larger than 6.5x can be expected for bigger problems.

\subsection{Weak scaling with constant \texorpdfstring{\(\delta\)}{delta}}\label{sec:weak-scaling-delta}

In a second set of experiments, we fix the interaction length \(\delta\) and vary the problem size \(L\).
This implies that the number of entries per matrix row will increase as \(L\) increases.
In order to keep the memory cost of the subdomain matrices \(A_{i}\) under control, this means that the number of subdomains per direction needs to scale like \(p \sim L^{2}\).

\begin{figure}
  \centering
  \includegraphics{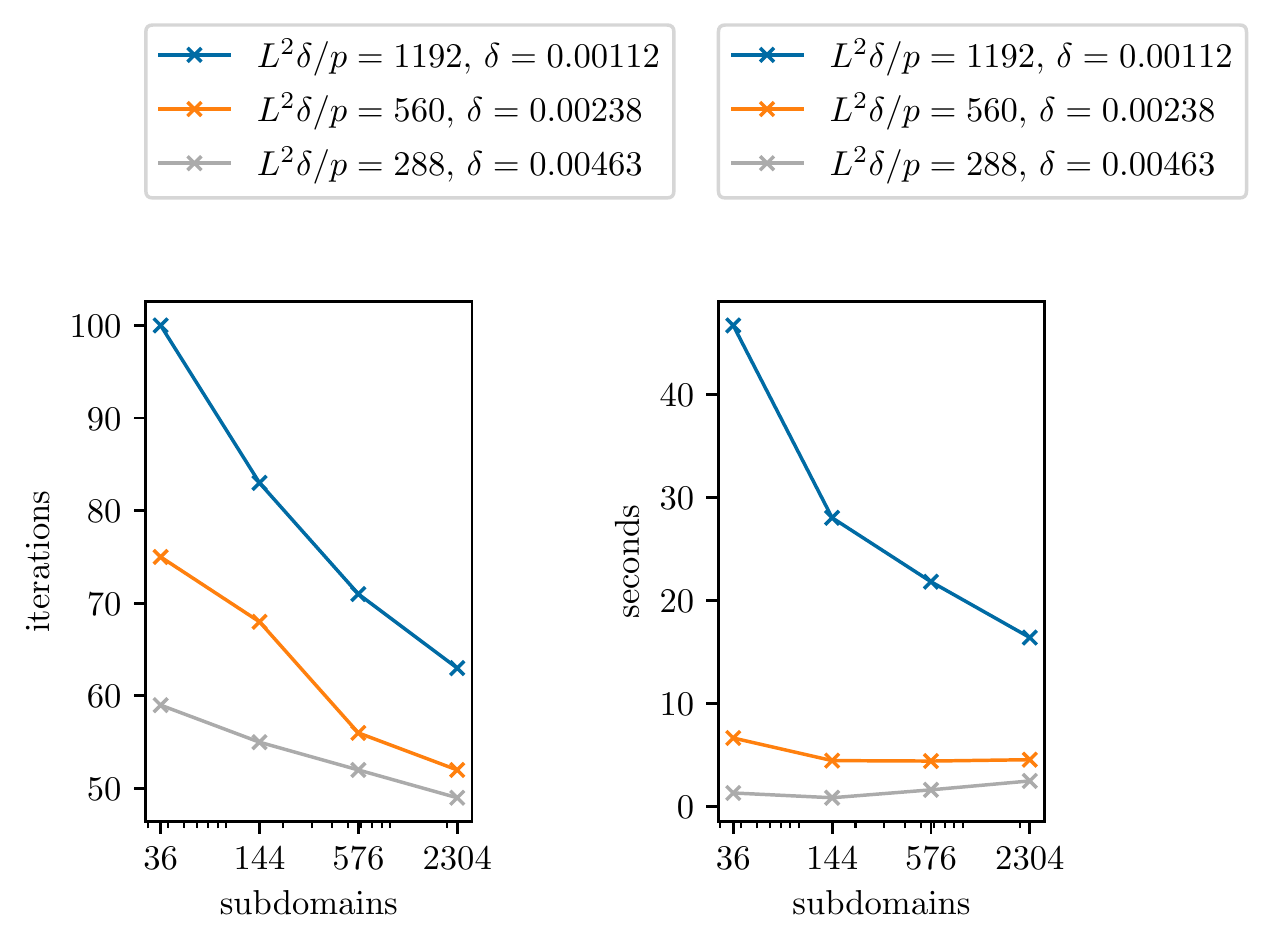}
  \caption{Weak scaling of the substructuring solver with constant \(\delta\) and \(p\sim L^{2}\).}
  \label{fig:weakScaling2}
\end{figure}
Iteration counts and timings results for the DD method with CG-Dirichlet preconditioner are shown in Figure~\ref{fig:weakScaling2}.
Since \(L/p\sim 1/L\), we observe decreasing iteration counts and solve times, as predicted by \eqref{eq:timeDD}.

\subsection{Strong scaling}
\label{sec:strong-scaling}

\begin{figure}
  \centering
  \includegraphics{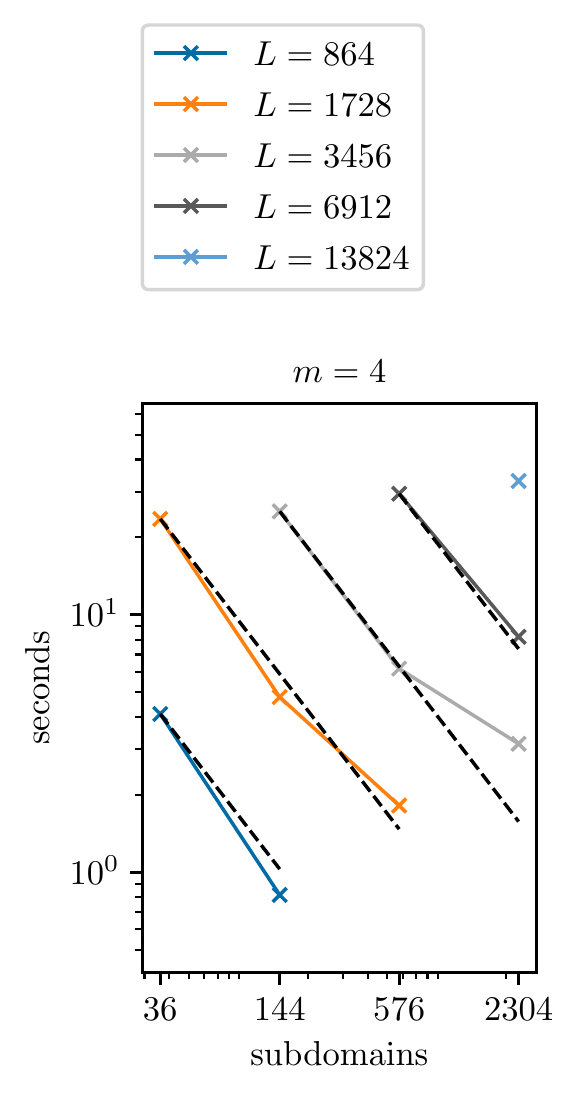}
  \includegraphics{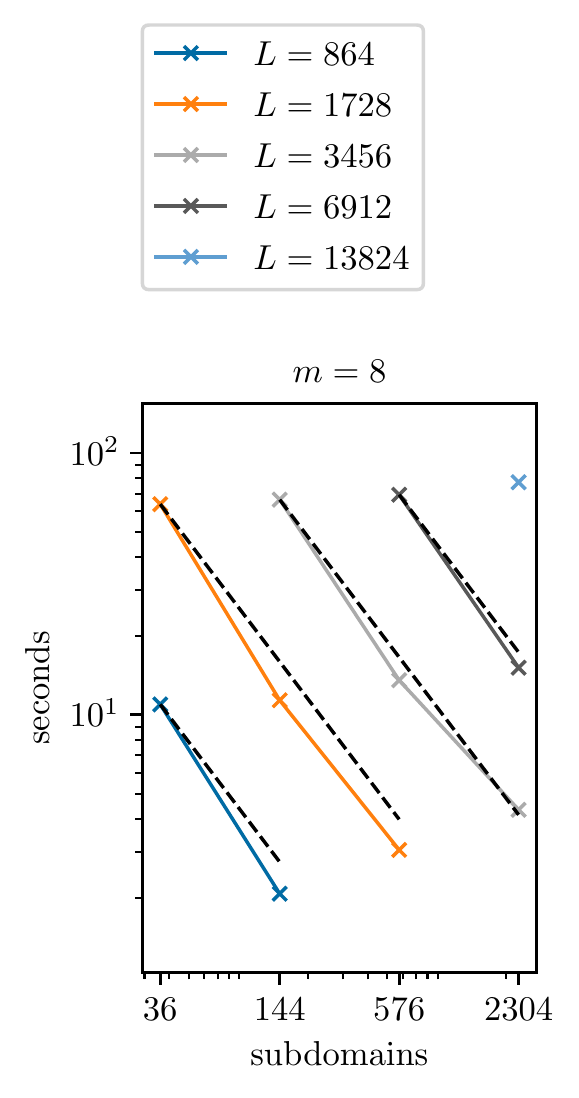}
  \caption{
    Strong scaling of the substructuring solver with constant \(L\), \(m\) and \(\delta=m/L\).
    The dashed lines are the ideal strong scaling.
  }
  \label{fig:strongScaling}
\end{figure}

In a final set of experiments, we explore the strong scaling properties of the substructuring solver.
For a fixed choice of \(L\), \(m\) and \(\delta=m/L\), we vary the number of subdomains per direction \(p\).
The resulting timings for \(m=4\) and \(m=8\) are shown in Figure~\ref{fig:strongScaling}.
As observed, the strong scaling is almost optimal.
For \(m=4\), we see the strong scaling performance taper off as the number of subdomains is increased by 16.
This can be explained by the fact that while computation is indeed divided by a factor of 16, communication is not and will start to dominate in the strong scaling limit.
Naturally, for \(m=8\) the amount of computation is higher, which explains the higher solve times compared to \(m=4\), but also the better strong scaling behavior.

\newpage

\section{Conclusion}\label{sec:conclusion}
We have presented a domain decomposition technique for the efficient numerical solution of nonlocal diffusion problems. 
Our method, based on a meshfree discretization of the variational multi-domain formulation that was previously introduced in \cite{Capodaglio2021}, solves the discretized problem in a FETI fashion and exploits effective preconditioning techniques introduced for local domain-decomposition problems. 

Our rigorous numerical studies highlight the fact that the proposed approach outperforms commonly used parallel solvers.
In fact, based on the scaling estimates derived in Section \ref{sec:numer-exper} and the corresponding numerical studies, we expect speedups (with respect to standard single-domain solvers) larger than 6.5x on problems of increasing size.  Furthermore, our tests in Sections \ref{sec:weak-scaling}--\ref{sec:strong-scaling} show excellent scalability properties that confirm our theoretical estimates. 

This work is the first rigorous numerical study in a two-dimensional multi-domain setting for nonlocal operators with finite horizon and, as such, it represents a fundamental step towards increasing the usability of nonlocal models in large scale simulations.  It is possible to generalize the proposed technique to any linear nonlocal problem, including nonlocal elasticity, in a straightforward way just by adapting the discretization matrices associated with the sub-problems.  Furthermore, due to the similarities between discretized nonlocal problems and systems resulting from long-stencil meshfree discretizations of PDEs, the impact of this work extends beyond nonlocal problems, as our algorithm can be applied without modifications to meshfree discretizations of PDEs.

\section*{Acknowledgment}

This work was supported by Sandia National Laboratories (SNL) Laboratory-directed Research and Development (LDRD) program, project 218318 and by the U.S. Department of Energy, Office of Advanced Scientific Computing Research under the Collaboratory on Mathematics and Physics-Informed Learning Machines for Multiscale and Multiphysics Problems (PhILMs) project. Sandia National Laboratories is a multimission laboratory managed and operated by National Technology and Engineering Solutions of Sandia, LLC., a wholly owned subsidiary of Honeywell International, Inc., for the U.S. Department of Energy's National Nuclear Security Administration contract number DE-NA0003525. This paper, SAND2021-5958 O, describes objective technical results and analysis. Any subjective views or opinions that might be expressed in the paper do not necessarily represent the views of the U.S. Department of Energy or the United States Government.

\bibliographystyle{abbrv}
\bibliography{mybib}
\end{document}